\documentclass{article}

\usepackage{proof}
\usepackage{latexsym}
\usepackage{amsfonts}
\usepackage{amssymb}

\newcommand{\brem}{\begin{remark}}
\newcommand{\erem}{\end{remark}}
\newcommand{\blem}{\begin{lemma}}
\newcommand{\elem}{\end{lemma}}
\newcommand{\bth}{\begin{theorem}}
\newcommand{\ethm}{\end{theorem}}
\newcommand{\benu}{\begin{enumerate}}
\newcommand{\eenu}{\end{enumerate}}
\newcommand{\bdes}{\begin{description}}
\newcommand{\edes}{\end{description}}
\newcommand{\bdf}{\begin{definition}}
\newcommand{\edf}{\end{definition}}
\newcommand{\bcor}{\begin{cor}}
\newcommand{\ecor}{\end{cor}}
\newcommand{\bprp}{\begin{proposition}}
\newcommand{\eprp}{\end{proposition}}
\newcommand{\bmlem}{\begin{mlemma}}
\newcommand{\emlem}{\end{mlemma}}
\newcommand{\bclm}{\begin{claim}}
\newcommand{\eclm}{\end{claim}}
\newcommand{\bprf}{{\bf Proof}.\hspace{2mm}}

\newcommand{\eprf}{\hspace*{\fill} $\Box$}

\newcommand{\beqn}{\begin{equation}}
\newcommand{\eeqn}{\end{equation}}
\newcommand{\beqnarr}{\begin{eqnarray}}
\newcommand{\eeqnarr}{\end{eqnarray}}
\newcommand{\beqnarrs}{\begin{eqnarray*}}
\newcommand{\eeqnarrs}{\end{eqnarray*}}

\newcommand{\spand}{\,\&\,}

\newtheorem{theorem}{Theorem}[section]
\newtheorem{definition}[theorem]{Definition}
\newtheorem{proposition}[theorem]{Proposition}
\newtheorem{lemma}[theorem]{Lemma}
\newtheorem{cor}[theorem]{Corollary}
\newtheorem{remark}[theorem]{Remark}
\newtheorem{mlemma}[theorem]{Main Lemma}
\newtheorem{claim}[theorem]{Claim}

\newcommand{\alp}{\alpha}

\newcommand{\veps}{\varepsilon}
\newcommand{\del}{\delta}
\newcommand{\Del}{\Delta}
\newcommand{\ome}{\omega}

\newcommand{\bet}{\beta}
\newcommand{\gam}{\gamma}
\newcommand{\Gam}{\Gamma}

\newcommand{\sig}{\sigma}
\newcommand{\Sig}{\Sigma}
\newcommand{\tht}{\theta}

\newcommand{\Lam}{\Lambda}
\newcommand{\vphi}{\varphi}

\newcommand{\fal}{\forall}
\newcommand{\exi}{\exists}

\newcommand{\Rarw }{\Rightarrow}

\newcommand{\Lrarw}{\Leftrightarrow}

\newcommand{\la}{\langle}
\newcommand{\ra}{\rangle}

\title{
Mahlo classes for first-order reflections
}

\author{Toshiyasu Arai
\\
Graduate School of Mathematical Sciences
\\
The University of Tokyo
\\
3-8-1 Komaba, Meguro-ku,
Tokyo 153-8914, JAPAN
\\
tosarai@ms.u-tokyo.ac.jp}
\date{}
\begin{document}

\maketitle

\begin{abstract}
In this note we axiomatize the $\Pi_{k+1}$-consequences in the set theory 
${\sf KP}\Pi_{N}$ for $\Pi_{N}$-reflecting universes
in terms of iterations of $\Pi_{i}$-recursively Mahlo operations
for $1\leq k\leq i<N$.
\end{abstract}

\section{Introduction}
For set-theoretic formulas $\varphi$, let
$
P\models\varphi:\Leftrightarrow (P,\in)\models\varphi
$, where
parameters occurring in $\vphi$ is understood to be in $P$.
In what follows, $V$ denotes a transitive and wellfounded model of 
${\sf KP}\ell$, which is a universe in discourse.
${\sf KP}\ell$ introduced in \cite{Jaeger86} denotes
a set theory for universes in which there
are unboundedly many admissible sets.
$ON$ denotes the class of ordinals in $V$.
$P, Q,\ldots$ denote non-empty transitive sets in $V\cup\{V\}$.

For $k=1,2,3,\ldots$, \textit{$\Pi_{k}$-recursively Mahlo operation}
 $M_{k}$ is defined as follows.
 Let $X\subset V$ be a class.
\beqnarrs
P\in M_{k}(X) & \Lrarw & 
P \mbox{ is } \Pi_{k}\mbox{-reflecting on } X\subset V
\\
& \Lrarw &
\mbox{for any }\Pi_{k}\mbox{-formula } \vphi
\mbox{ with parameters from } P,
\\
&&
\mbox{ if } P\models\vphi,
\mbox{ then }
Q\models\vphi
\mbox{ for a } Q\in X\cap P
\eeqnarrs
When $P\in M_{k}:\Lrarw 
P\in M_{k}(V)$, $P$ is said to be $\Pi_{k}$-reflecting universe.
Note that if $P$ is $\Pi_{k}$-reflecting on $X$, then
it is $\Sig_{k+1}$-reflecting on $X$.

For ordinals $\pi\in ON$, $\pi$ is said to be $\Pi_{k}$-reflecting on $X\subset ON$
iff $L_{\pi}\in M_{k}(\{L_{\alp}:\alp\in X\})$.

It is easy to see from the existence of a universal $\Pi_{k}$-formula
that $P\in M_{k}(X)$ is a $\Pi_{k+1}(X)$-formula
uniformly on admissible sets $P$:
there exists an $X$-positive $\Pi_{k+1}(X^{+})$-formula 
$\tht_{k+1}(X)$ such that
for any admissible set $P$,
$P\in M_{k}(\mathcal{X})$ iff $(P,\mathcal{X}\cap P)\models \tht_{k+1}(X)$.

An iteration of the operation $X\mapsto M_{k}(X)$ 
along a $\Del$-relation $\prec$ is defined by:
\[
P\in M_{k}(a;\prec)\Lrarw P\in\bigcap\{M_{k}(M_{k}(b;\prec)): P\ni b\prec a\}
.\]
Again we see that there exists a $\Pi_{k+1}$-formula $\sig_{k+1}(a)$ such that for any admissible set $P$ and $a\in P$,
$P\in M_{k}(a;\prec)$ iff $P\models \sig_{k+1}(a)$.

Let $ON^{\veps}\subset V$ and $<^{\veps}$ be $\Del$-predicates such that
for any transitive and wellfounded model $V$ of 
the Kripke-Platek set theory $\mbox{{\sf KP}}\ome$ with the axiom of Infinity,
$<^{\veps}$ is a well ordering of type $\veps_{\mathbb{K}+1}$
 on $ON^{\veps}$
for the order type $\mathbb{K}$ of the class $ON$ in $V$.
$ON^{\veps}$ is the class of codes $\lceil\alp\rceil$ of ordinals 
$\alp<\veps_{\mathbb{K}+1}$.
Each admissible set $P$ is assumed to be closed under addition
$(\lceil\alp\rceil,\lceil\bet\rceil)\mapsto\lceil\alp+\bet\rceil$ and
exponentiation $\lceil\alp\rceil\mapsto\lceil\ome^{\alp}\rceil$ on codes
besides $\{\lceil\mathbb{K}\rceil\}\cup\{\lceil\alp\rceil: \alp\in P\cap ON\}$.
In other words, the order type of the set of codes $\lceil\alp\rceil$ in $P$
is the next epsilon number above
the order type of $P\cap ON$.

$<^{\veps}$ is assumed to be a canonical ordering such that 
$\mbox{{\sf KP}}\ome$ proves the fact that $<^{\veps}$ is a linear ordering, and for any formula $\vphi$
and each $n<\ome$,
\beqn\label{eq:trindveps}
\mbox{{\sf KP}}\ome\vdash\fal x(\fal y<^{\veps}x\,\vphi(y)\to\vphi(x)) \to \fal x<^{\veps}\lceil\ome_{n}(\mathbb{K}+1)\rceil\vphi(x)
\eeqn
In what follows let us write $a<b$ for $a<^{\veps}b$ and $a,b\in ON^{\veps}$
when no confusion likely occurs.

By induction on ordinals $a<\veps_{\mathbb{K}+1}$(the next epsilon number above $\mathbb{K}$) we see that
$M_{k+1}
\subset \bigcap_{a<\veps_{\mathbb{K}+1}}M_{k}(a;<^{\veps})
$.
In \cite{consv} it is shown that
a set theory ${\sf KP}\Pi_{N}$ for $\Pi_{N}$-reflecting universes
is $\Pi_{N}$-conservative over the theory
${\sf KP}\ome+\{V\in M_{N-1}(a; <^{\veps}): a\in ON^{\veps}\}$.
In other words, the set $\{V\in M_{N-1}(a; <^{\veps}): a\in ON^{\veps}\}$ axiomatizes
the $\Pi_{N}$-consequences of ${\sf KP}\Pi_{N}$ over ${\sf KP}\ome$.

In this note we axiomatize the $\Pi_{k+1}$-consequences of ${\sf KP}\Pi_{N}$
for $1\leq k< N$ through $\Pi_{k+1}$-formulas $V\in M_{k}(\vec{a})$
for a ramified iteration of operations $M_{i}$.

\section{Ramified iterations of recursively Mahlo operations}
Inspired from M. Rathjen\cite{Rathjen94} 
let us introduce a ramified iterations of $M_{k}$.
Let $N\geq 2$.
In what follows let $\Lam=\veps_{\mathbb{K}+1}$, and
$\alp,\bet,\gam$ denote ordinals below $\Lam$.

\bdf\label{df:uppertriangle}
{\rm
Let $\vec{\alp}=\la\alp_{k},\ldots,\alp_{N-1}\ra$ and $\vec{\bet}=\la\bet_{k},\ldots,\bet_{N-1}\ra$
be sequences of ordinals $\alp_{i},\bet_{i}<\Lam$ in the same length 
$lh(\vec{\alp})=lh(\vec{\bet})=N-k>0$.
\benu
\item
$\vec{\bet}<\vec{\alp}:\Lrarw \fal i(\bet_{i}<\alp_{i})$.
\item
$\vec{\bet}\bullet\vec{\alp}=\la\vec{\gam}_{k},\ldots,\vec{\gam}_{N-1}\ra$ denotes the sequence of 
sequences $\vec{\gam}_{i}$ of ordinals 
in length $lh(\vec{\gam}_{i})=N-i$
defined by
$\vec{\gam}_{i}=\la\bet_{i},\alp_{i+1},\ldots,\alp_{N-1}\ra$.
\eenu
}
\edf
Each sequence $\vec{\gam}_{i}$ in the sequence 
$\vec{\bet}\bullet\vec{\alp}$ is the $(i-k+1)$-th row
in the upper triangular matrix:
\[
\left(
\begin{array}{c}
\vec{\gam}_{k}
\\
\vec{\gam}_{k+1}
\\
\vdots
\\
\vec{\gam}_{N-2}
\\
\vec{\gam}_{N-1}
\end{array}
\right)
=
\left(
\begin{array}{ccccc}
\bet_{k} & \alp_{k+1} & \alp_{k+2} & \cdots & \alp_{N-1}
\\
 & \bet_{k+1} & \alp_{k+2} & \cdots & \alp_{N-1}
 \\
 & & \ddots & & \vdots
 \\
 & & & \bet_{N-2} & \alp_{N-1}
 \\
 & & & & \bet_{N-1}
\end{array}
\right)
\]

\bdf\label{df:Mhclass}
{\rm
For $1\leq k<N$ and sequences $\vec{\alp}=\la\alp_{k},\ldots,\alp_{N-1}\ra$ of ordinals $\alp_{i}<\Lam$ in length $N-k$, let us define a class $Mh_{k}(\vec{\alp})$ as follows.
\[
P\in 
Mh_{k}(\vec{\alp}) :\Lrarw \vec{\alp}\in P \spand 
\fal \vec{\bet}\in P
\left[
\vec{\bet}<\vec{\alp} \Rarw P\in M_{k}\left(Mh_{k}(\vec{\bet}\bullet\vec{\alp})\right)
\right]
\]
where 
$\vec{\alp}\in P :\Lrarw
\fal i(\alp_{i}\in P)$, and
$Mh_{k}(\la\vec{\gam}_{k},\ldots,\vec{\gam}_{N-1}\ra):=\bigcap_{i\geq k} Mh_{i}(\vec{\gam}_{i})$.

When $k=N$, i.e., $\vec{\alp}$ is the empty sequence $\emptyset$,
set
$Mh_{N}(\emptyset):=M_{N}$.
}
\edf
For example $Mh_{N-1}(\la\alp\ra)=M_{N-1}(\alp)$ for ordinals $\alp<\Lam$, and
for $\alp_{N-2},\alp_{N-1}\in P$,
$P\in Mh_{N-2}(\la\alp_{N-2},\alp_{N-1}\ra)$ iff 
for any $\bet_{N-2}\in P\cap \alp_{N-2}$
and any $\bet_{N-1}\in P\cap \alp_{N-1}$,
$P$ is $\Pi_{N-2}$-reflecting on the intersection
$Mh_{N-2}(\la\bet_{N-2},\alp_{N-1}\ra)\cap M_{N-1}(\bet_{N-1})$, i.e.,
$P\in M_{N-2}
\left(
Mh_{N-2}(\la\bet_{N-2},\alp_{N-1}\ra)\cap M_{N-1}(\bet_{N-1})
\right)$ holds.

Again we see that there exists a $\Pi_{k+1}$-formula $\sig_{k+1}(\vec{a})$ such that for any admissible set $P$ and $\vec{\alp}\in P$,
$P\in Mh_{k}(\vec{\alp})$ iff $P\models \sig_{k+1}(\vec{\alp})$.

\bdf\label{df:TNk}
{\rm
For $1\leq k\leq N$, let $T^{(N)}_{k}$ denote an extension of
${\sf KP}\ome$ by the axioms $V\in Mh_{k}(\vec{\alp})$ with
$\vec{\alp}=\la\alp_{k},\ldots,\alp_{N-1}\ra$ for $\alp_{i}<\Lam$.
\[
T^{(N)}_{k}:= {\sf KP}\ome+\{V\in Mh_{k}(\vec{\alp}): lh(\vec{\alp})=N-k,
\fal i\geq k(\alp_{i}<\Lam)\}
.\]
}
\edf
$T^{(N)}_{N}$ is defined to be 
the set theory ${\sf KP}\Pi_{N}$ for $\Pi_{N}$-reflecting universes
since $Mh_{N}(\emptyset)=M_{N}$.

We show the following Theorem \ref{th:cnsvlstep}.

\bth\label{th:cnsvlstep}
For each $1\leq k<N$, $T^{(N)}_{k+1}$ is $\Pi_{k+1}$-conservative over 
$T^{(N)}_{k}$.
\end{theorem}

\bcor\label{cor:consvlstep}
$T^{(N)}_{N}={\sf KP}\Pi_{N}$ is $\Pi_{3}$-conservative over $T^{(N)}_{2}$.
Moreover
$T^{(N)}_{N}$ is $\Pi_{2}$-conservative over $T^{(N)}_{1}$.
\ecor

First let us show that $T^{(N)}_{k}$ is a subtheory of $T^{(N)}_{k+1}$.

\blem\label{lem:prov}
$T^{(N)}_{k+1}$ proves $V\in Mh_{k}(\vec{\alp})$ for {\rm each}
$\vec{\alp}$ with $lh(\vec{\alp})=N-k$.
\elem
\bprf
We show that $T^{(N)}_{k+1}$ proves 
$\fal\vec{\alp}<\ome_{n}(\mathbb{K}+1)(V\in Mh_{k}(\vec{\alp}))$ for each $n<\ome$,
where
$\la\alp_{k},\ldots,\alp_{N-1}\ra<\ome_{n}(\mathbb{K}+1):\Lrarw 
\fal i\geq k(\alp_{i}<\ome_{n}(\mathbb{K}+1))$.
Let $\vec{\alp}=\la\alp_{k+1},\ldots,\alp_{N-1}\ra$ be a sequence of ordinals in length $N-k-1$.
Argue in $T^{(N)}_{k+1}$.
We have $V\in Mh_{k+1}(\vec{\alp})$.
We show by induction on ordinals $\alp_{k}<\ome_{n}(\mathbb{K}+1)$,
cf.\,(\ref{eq:trindveps}) that
$\fal \alp_{k}<\ome_{n}\left(\mathbb{K}+1)(V\in Mh_{k}(\la\alp_{k}\ra*\vec{\alp})\right)$,
where $\la\alp_{k}\ra*\vec{\alp}=(\alp_{k},\alp_{k+1},\ldots,\alp_{N-1})$.

Let $\la\bet_{k}\ra*\vec{\bet}<\la\alp_{k}\ra*\vec{\alp}$.
We obtain $\vec{\bet}<\vec{\alp}$ and $V\in M_{k+1}(Mh_{k+1}(\vec{\bet}\bullet\vec{\alp}))$
by $V\in Mh_{k+1}(\vec{\alp})$.
On the other hand we have $V\in Mh_{k}(\la\bet_{k}\ra*\vec{\alp})$ by IH and $\bet_{k}<\alp_{k}$.
Since this is a $\Pi_{k+1}$-sentence holding on 
$V$, which is in $M_{k+1}(Mh_{k+1}(\vec{\bet}\bullet\vec{\alp}))$, we obtain
$V\in M_{k+1}(Mh_{k}(\la\bet_{k}\ra*\vec{\alp})\cap Mh_{k+1}(\vec{\bet}\bullet\vec{\alp}))$.
On the other side
we see that 
$(\la\bet_{k}\ra*\vec{\bet})\bullet(\la\alp_{k}\ra*\vec{\alp})=
\la \la\bet_{k}\ra*\vec{\alp}\ra*(\vec{\bet}\bullet\vec{\alp})$.
Therefore
$Mh_{k}(\la\bet_{k}\ra*\vec{\alp})\cap Mh_{k+1}(\vec{\bet}\bullet\vec{\alp})
= Mh_{k}\left( (\la\bet_{k}\ra*\vec{\bet})\bullet(\la\alp_{k}\ra*\vec{\alp})\right)$.
\eprf

\section{Conservativity}

Assume that $T^{(N)}_{k+1}\vdash A$ for a $\Pi_{k+1}$-sentence $A$.
Pick a sequence $\vec{\alp}$ of ordinals in length $N-k-1$
such that
${\sf KP}\ome\vdash V\in Mh_{k+1}(\vec{\alp}) \to A$.
We need to find a longer sequence $\bet$ such that
${\sf KP}\ome\vdash V\in Mh_{k}(\vec{\bet}) \to A$.

In what follows we work in the intuitionistic fixed point theory 
${\rm FiX}^{i}(T^{(N)}_{k})$ over $T^{(N)}_{k}$.
In \cite{Arai15} it is shown that
${\rm FiX}^{i}(T^{(N)}_{k})$ is a conservative extension of $T^{(N)}_{k}$.

As in \cite{consv}
${\sf KP}\ome+(V\in Mh_{k+1}(\vec{\alp}))$ is embedded to an infinitary 
one-sided sequent calculus 
with inference rules 
$(Mh_{k+1}(\vec{\alp},\vec{\bet}))$.
In one-sided sequent calculi, formulas are generated from atomic formulas and their negations
$a\in b,a\not\in b$
by propositional connectives $\lor,\land$ and quantifiers $\exi,\fal$.
It is convenient for us to have bounded quantifications $\exi x\in a,\fal x\in a$ besides unbounded ones $\exi x,\fal x$.
The negation $\lnot A$ of formulas $A$ is defined recursively by de Morgan's law and elimination of double negations.
Also $(A\to B):\equiv(\lnot A\lor B)$.
The language of the calculus is
 a set-theoretic language $\mathcal{L}_{V}$ with names
$c_{a}$ for each set $a\in V$.
Let us identify the name $c_{a}$ with the set $a$.

A finite set $\Gamma$ of sentences, a \textit{sequent}
 in the language $\mathcal{L}_{V}$
is intended to denote the disjunction $\bigvee\Gamma:=\bigvee\{A:A\in\Gamma\}$.
A sequent $\Gamma$ is said to be \textit{true} in $P\in V\cup\{V\}$ 
iff $\bigvee\Gamma$ is true in $P$.

We assign disjunctions or conjunctions to sentences as follows.

\bdf\label{df:assigndc}
{\rm
\begin{enumerate}
\item
For a $\Delta_{0}$-sentence $M$
\[
M:\simeq
\left\{
\begin{array}{ll}
\bigvee(A_{\iota})_{\iota\in J} & \mbox{{\rm if }} M \mbox{ {\rm is false in }} V
\\
\bigwedge(A_{\iota})_{\iota\in J} &  \mbox{{\rm if }} M \mbox{ {\rm is true in }} V
\end{array}
\right.
\mbox{{\rm with }} J:=\emptyset
.\]

In what follows we consider the unbounded sentences.

\item
$(A_{0}\lor A_{1}):\simeq\bigvee(A_{\iota})_{\iota\in J}$
and
$(A_{0}\land A_{1}):\simeq\bigwedge(A_{\iota})_{\iota\in J}$
with $J:=2$.

\item
$
\exists x\in a\, A(x):\simeq\bigvee(A(b))_{b\in J}$
and
$
\forall x\in a\, A(x):\simeq\bigwedge(A(b))_{b\in J}
$ 
with
$
J:=a$.

\item
$
\exists x\, A(x):\simeq\bigvee(A(b))_{b\in J}
$ 
and
$
\forall x\, A(x):\simeq\bigwedge(A(b))_{b\in J}
$ 
with
$
J:=V
$.

\end{enumerate}
}
\edf

\bdf\label{df:dp}
{\rm 
The \textit{depth} $\mbox{{\rm dp}}(A)<\omega$ of
 $\mathcal{L}_{V}$-sentences $A$ is defined recursively as follows.
\begin{enumerate}
\item
$\mbox{{\rm dp}}(A)=0$ if $A\in\Delta_{0}$.

In what follows we consider unbounded sentences $A$.

\item
$\mbox{{\rm dp}}(A)=\max\{\mbox{{\rm dp}}(A_{i}):i<2\}+1$ if
 $A\equiv (A_{0}\circ A_{1})$ 
for $\circ\in\{\lor,\land\}$.

\item
$\mbox{{\rm dp}}(A)=\mbox{{\rm dp}}(B(\emptyset))+1$ 
if $A\in\{ (Q x \, B(x)), (Q x\in a\, B(x)) :a\in V\}$
for $Q\in\{\exists,\forall\}$.
\end{enumerate}
}
\edf

\bdf\label{df:sfk}
{\rm
\begin{enumerate}
\item
For $\mathcal{L}_{V}$-sentences $A$,
$
{\sf k}(A):=\{a\in V: c_{a}\mbox{ {\rm occurs in} } A\}
$.

\item
For sets $\Gamma$ of sentences,
 ${\sf k}(\Gamma):=\bigcup\{{\sf k}(A):A\in\Gamma\}$.

\item
For $\iota\in V$ and a transitive model $P\in V$
of ${\sf KP}\omega$,
$P(\iota)\in V\cup\{V\}$ 
denotes the smallest transitive model of
${\sf KP}\omega$ such that $P\cup\{\iota\}\subset P(\iota)$.
Note that $V$ is assumed to be a model of ${\sf KP}\ell$.
\end{enumerate}
}
\edf

Inspired by operator controlled derivations due to W. Buchholz \cite{Buchholz},
let us define a relation $P\vdash^{\alpha}_{c}\Gamma$ for transitive models $P\in V\cup\{V\}$ of ${\sf KP}\omega$.
The relation $P\vdash^{\alpha}_{c}\Gamma$ is defined as a fixed point of a strictly positive formula $H$
\[
H(P,\alpha,c,\Gamma)\Leftrightarrow P\vdash^{\alpha}_{c}\Gamma
\]
in $\mbox{FiX}^{i}({\sf KP}\ell)$.

Note that $P$ 
is closed under $a\mapsto rank(a)$ for $rank(a)=\sup\{rank(b)+1:b\in a\}$.

\bdf\label{df:controlderreg}
{\rm 
Fix a sequence $\vec{\alp}\neq\la 0,\ldots,0\ra$ of ordinals in length $N-k-1$.
Let $P\in V\cup\{V\}$ be a transitive model of ${\sf KP}\omega$,  codes
$\alp<\Lam$ and $c<\omega$.

$P\vdash^{\alp}_{c}\Gamma$ holds if
\begin{equation}\label{eq:controlord}
{\sf k}(\Gamma)\cup\{\alp\}\subset P
\end{equation}
and one of the following
cases holds:

\begin{description}
\item[$(\bigvee)$]
There is an $A\in\Gamma$ such that
$A\simeq\bigvee(A_{\iota})_{\iota\in J}$, and for an $\iota\in J$  and an
 $\alp(\iota)<\alp$,
$P\vdash^{\alp(\iota)}_{c}\Gamma,A_{\iota}$.
\[
\infer[(\bigvee)]
{P\vdash^{\alp}_{c}\Gamma}{P\vdash^{\alp(\iota)}_{c}\Gamma,A_{\iota}}
\]

\item[$(\bigwedge)$]
There is an $A\in\Gamma$ such that
$A\simeq\bigwedge(A_{\iota})_{\iota\in J}$, and for any $\iota\in J$, 
there is an $\alp(\iota)$ such that $\alp(\iota)<\alp$ and
$P(\iota)\vdash^{\alp(\iota)}_{c}\Gamma,A_{\iota}$.
\[
\infer[(\bigwedge)]{P\vdash^{\alp}_{c}\Gamma}{\{P(\iota)\vdash^{\alp(\iota)}_{c}\Gamma,A_{\iota}:\iota\in J\}}
\]

\item[$(cut)$]
There are $C$ and $\alp_{0}<\alp$ such that
$\mbox{{\rm dp}}(C)<c$, and
$P\vdash^{\alp_{0}}_{c}\Gamma,\lnot C$ and
 $P\vdash^{\alp_{0}}_{c}C,\Gamma$.

\[
\infer[(cut)]{P\vdash^{\alp}_{c}\Gamma}
{P\vdash^{\alp_{0}}_{c}\Gamma,\lnot C & P\vdash^{\alp_{0}}_{c}C,\Gamma}
\]

\item[$\left(Mh_{k+1}(\vec{\alp},\vec{\bet})\right)$]
There are a sequence $\vec{\bet}$ of ordinals in length $N-k-1$,
a sequent $\Del\subset\Sig_{k+1}$ and
$\alp_{\ell},\alp_{r}<\alp$ for which the followings hold:
$\vec{\bet}<\vec{a}$ and $\vec{\bet}\cup\vec{\alp}\subset P$.
$P\vdash^{\alp_{\ell}}_{c}\Gamma,\lnot\del$ for each $\del\in\Del$.

For each $Q\in Mh_{k+1}(\vec{\bet}\bullet\vec{\alp})$,
$P(Q)\vdash^{\alp_{r}}_{c}\Gamma,\Del^{(Q)}$ holds, where
for sentences $\del$, $\del^{(Q)}$ denotes the result of
restricting each unbounded quantifier $\exi x,\fal x$ in $\del$
to $\exi x\in Q,\fal x\in Q$, resp.
$\Del^{(Q)}=\{\del^{(Q)}:\del\in\Del\}$.
\[
\infer[\left(Mh_{k+1}(\vec{\alp},\vec{\bet})\right)]{P\vdash^{\alp}_{c}\Gamma}
{
P\vdash^{\alp_{\ell}}_{c}\Gamma,\lnot\del \, (\del\in\Del)
&
P(Q)\vdash^{\alp_{r}}_{c}\Gamma,\Del^{(Q)}
\,(Q\in Mh_{k+1}(\vec{\bet}\bullet\vec{\alp}))
}
\]

\end{description}
}
\edf

\blem\label{th:embedreg}{\rm (Embedding)}\\
If ${\sf KP}\ome+(V\in Mh_{k+1}(\vec{\alp})) \vdash A$ for sentences $A$, then
there exist $p,c<\ome$ such that
 $P\vdash_{c}^{\mathbb{K}\cdot p}A$ holds 
 for any transitive model $P$ of {\sf KP}$\ome$.
\elem
\bprf
As in \cite{consv} we see by induction on $rank(a)<\mathbb{K}$ that 
$P(a)\vdash_{0}^{2d+3 rank(a)}B, \forall x\in a\, A(x)$
for $d=\mbox{{\rm dp}}(A)$ and $B\equiv(\lnot\forall x(\forall y\in x\, A(y)\to A(x)))$.

Next consider the axioms in ${\sf KP}\ome$ other than Foundation and
$\Del_{0}$-Collection.
It is a $\Pi_{2}$-axiom $\forall x,y\exists z\, \tht(x,y,z)$.
Let $a,b\in V$.
Since $P(a,b)$ is a transitive model of ${\sf KP}\omega$ and $a,b\in P(a,b)$,
pick a $c\in P(a,b)$ such that the $\Delta_{0}$-formula $\tht(a,b,c)$ holds in 
$P(a,b)$, and in $V$.
Since this is a true $\Delta_{0}$-sentence, we have
 $P(a,b)\vdash^{0}_{0}\tht(a,b,c)$, and
$P\vdash^{3}_{0} \forall x,y\exists z\, \tht(x,y,z)$.

To get a derivation of $V\in Mh_{k+1}(\vec{\alp})$, let
$\vec{\bet}<\vec{\alp}$,
and $\del\in\Sig_{k+1}$ with ${\sf k}(\del)\subset P$.
We obtain $P\vdash^{2d}_{0}\del,\lnot\del$ for $d=\mbox{dg}(\del)$ by tautology.
Let $Q\in Mh_{k+1}(\vec{\bet}\bullet\vec{\alp})\cap V$.
We obtain
$P(Q)\vdash^{0}_{0} \del^{(Q)}, \lnot\del^{(Q)}$.
Moreover
$P(Q)\vdash^{0}_{0}Q\in Mh_{k+1}(\vec{\bet}\bullet\vec{\alp})$
since $Q\in Mh_{k+1}(\vec{\bet}\bullet\vec{\alp})$ is a true $\Del_{0}$-sentence.
Hence
$P(Q)\vdash^{1}_{0} \del^{(Q)},\exi x\in Mh_{k+1}(\vec{\bet}\bullet\vec{\alp}) \lnot\del^{(x)}$.
An inference rule $(Mh_{k+1}(\vec{\alp},\vec{\bet}))$ yields
$P\vdash^{\ome}_{0} \del, \exi x\in Mh_{k+1}(\vec{\bet}\bullet\vec{\alp}) \lnot\del^{(x)}$.

$\Del_{0}$-Collection follows from $V\in Mh_{k+1}(\vec{\alp})$ for $k\geq 1$.
\eprf

\blem\label{lem:predcereg}{\rm (Predicative Cut-elimination)}\\
For {\rm each} $n<\ome$ the following holds:

For any transitive model $P$ of {\sf KP}$\ome$, $c<\ome$ and 
$\alp<\ome_{n}(\mathbb{K}+1)$,
$P\vdash^{\alp}_{c+1}\Gam
 \Rarw 
P\vdash^{\ome^{\alp}}_{c}\Gam$.
\elem
\bprf
This is seen from the fact that
if $P\vdash^{\alp}_{c}\Gam,\lnot C$ and $P\vdash^{\bet}_{c}C,\Del$
with $\mbox{dg}(C)\leq c$, then
$P\vdash^{\alp+\bet}_{c}\Gam,\Del$, where
$C$ is either a $\Del_{0}$-formula or
an existential formula.
The fact is shown by induction on $\bet<\ome_{n}(\mathbb{K}+1)$.
\eprf

\blem\label{lem:CollapsingthmKR}{\rm (Elimination of 
$(Mh_{k+1}(\vec{\alp},\vec{\bet}))$)}\\
For {\rm each} $n<\ome$ the following holds:

Let $\Gam\subset\Pi_{k+1}$.
Suppose $P_{0}\vdash^{\alp}_{0}\Gam$ with $\alp<\ome_{n}(\mathbb{K}+1)$
 and $P_{0}\in P\in Mh_{k}(\la\alp\ra*\vec{\alp})$.
Then $\Gam$ is true in $P$.
\elem
\bprf
We show the lemma by induction on $\alp<\ome_{n}(\mathbb{K}+1)$, 
cf.\,(\ref{eq:trindveps}).
Consider the case when the last inference is an $(Mh_{k+1}(\vec{\alp},\vec{\bet})))$
for $\vec{\bet}<\vec{\alp}$.
There are $\Del\subset\Sig_{k+1}$ and
$\alp_{\ell},\alp_{r}<\alp$ such that
\[
\infer[\left(Mh_{k+1}(\vec{\alp},\vec{\bet})\right)]{P_{0}\vdash^{\alp}_{0}\Gamma}
{
P_{0}\vdash^{\alp_{\ell}}_{0}\Gamma,\lnot\del \, (\del\in\Del)
&
P_{0}(Q)\vdash^{\alp_{r}}_{0}\Gamma,\Del^{(Q)}
\,(Q\in Mh_{k+1}(\vec{\bet}\bullet\vec{\alp}))
}
\]
We obtain $P\in Mh_{k}(\la\alp\ra*\vec{\alp})\subset Mh_{k}(\la\alp_{r}\ra*\vec{\alp})$
by $\alp_{r}<\alp$.
Assume
\beqn\label{eq:elim0}
P\not\models\bigvee\Gam
\eeqn
IH yields $P\models \bigvee\Del^{(Q)}$, i.e.,
\beqn\label{eq:elim1}
\fal Q\in  Mh_{k+1}(\vec{\bet}\bullet\vec{\alp})\cap P
\left[
\bigvee\Del^{(Q)}
\right]
\eeqn
On the other side let
$P_{0}\in Q\in Mh_{k}(\la\alp_{\ell}\ra*\vec{\alp})$.
By IH we obtain 
\beqn\label{eq:elim2}
\fal\del\in\Del
\left[
Q\models\bigvee\Gam \lor Q\models \lnot\del
\right]
\eeqn

Now let
$P_{0}\in Q\in Mh_{k}(\la\alp_{\ell}\ra*\vec{\alp})\cap Mh_{k+1}(\vec{\bet}\bullet\vec{\alp})\cap P$.
From (\ref{eq:elim1}) there exists a $\del\in\Del$ such that
$Q\models\del$.
For this $\del$ we obtain $Q\models\bigvee\Gam$ by (\ref{eq:elim2}).

Thus we have shown that 
\beqn\label{eq:elim3}
P_{0}\in Q\in Mh_{k}(\la\alp_{\ell}\ra*\vec{\alp})\cap Mh_{k+1}(\vec{\bet}\bullet\vec{\alp})\cap P
\Rarw
Q\models\bigvee\Gam
\eeqn

On the other side 
$P\in Mh_{k}(\la\alp\ra*\vec{\alp})\subset 
M_{k}\left(Mh_{k}(\la\alp_{\ell}\ra*\vec{\alp})\cap Mh_{k+1}(\vec{\bet}\bullet\vec{\alp})\right)$ 
holds
by $\alp_{\ell}<\alp$.
Namely we have $\la\alp_{\ell}\ra*\vec{\bet}<\la\alp\ra*\vec{\alp}$
and $(\la\alp_{\ell}\ra*\vec{\bet})\bullet(\la\alp\ra*\vec{\alp})=
\la\la\alp_{\ell}\ra*\vec{\alp}\ra*(\vec{\bet}\bullet\vec{\alp})$,
and hence
$Mh_{k}(\la\alp_{\ell}\ra*\vec{\alp})\cap Mh_{k+1}(\vec{\bet}\bullet\vec{\alp})=
Mh_{k}\left((\la\alp_{\ell}\ra*\vec{\bet})\bullet(\la\alp\ra*\vec{\alp})\right)$.

From the assumption (\ref{eq:elim0}), $P\models\lnot\bigvee\Gam$ for
$\lnot\bigvee\Gam\in\Sig_{k+1}$, we see that
 there exists a $Q$ such that
$P_{0}\in Q\in Mh_{k}((\alp_{\ell})*\vec{\alp})\cap Mh_{k+1}(\vec{\bet}\bullet\vec{\alp})\cap P$
and $Q\models\lnot\bigvee\Gam$.
This contradicts with (\ref{eq:elim3}).
Therefore the assumption (\ref{eq:elim0}) is refuted, and
we conclude $P\models\bigvee\Gam$.

Other cases are seen easily from IH as in \cite{consv}.
\eprf
\\

\noindent
{\bf Proof} of Theorem \ref{th:cnsvlstep}.
Argue in the intuitionistic fixed point theory 
${\rm FiX}^{i}(T^{(N)}_{k})$ over $T^{(N)}_{k}$.
Assume that $T^{(N)}_{k+1}\vdash A$ for a $\Pi_{k+1}$-sentence $A$.
Pick a non-zero sequence $\vec{\alp}$ such that
${\sf KP}\ome\vdash V\in Mh_{k+1}(\vec{\alp}) \to A$ holds.
Embedding \ref{th:embedreg} yields 
$L_{\ome_{1}^{CK}}\vdash^{\mathbb{K}\cdot p}_{c}A$ for some $p,c<\ome$.
Predicative Cut-elimination \ref{lem:predcereg} yields
$L_{\ome_{1}^{CK}}\vdash^{\alp}_{0}A$ with an $\alp<\ome_{c+1}(\mathbb{K}+1)$.
From Lemma \ref{lem:CollapsingthmKR}  we see that $A$ is true in $V$ such that
$L_{\ome_{1}^{CK}}\in V\in Mh_{k}((\alp)*\vec{\alp})$.
Therefore ${\rm FiX}^{i}(T^{(N)}_{k})\vdash A$.
We conclude $T^{(N)}_{k}\vdash A$ by \cite{Arai15}.
\\

\noindent
From Corollary \ref{cor:consvlstep} we see that
$T^{(N)}_{N}={\sf KP}\Pi_{N}$ is $\Pi_{3}$-conservative over $T^{(N)}_{2}$
[$\Pi_{2}$-conservative over $T^{(N)}_{1}$], resp.
However 
these reductions seem to be useless in
an ordinal analysis of the set theory ${\sf KP}\Pi_{N}$,
since both of the classes $Mh_{2}(\vec{\alp})$ and $Mh_{1}(\vec{\alp})$
involve the higher Mahlo operations $M_{k}$ with $k>2$
when $N>3$.
In the ordinal analysis of ${\sf KP}\Pi_{N}$ in \cite{KPpiNsmpl}, 
Mahlo classes $Mh_{k}^{a}(\xi)$ are introduced.
These classes are indexed by ordinals $\xi$ up to
the next epsilon number $\veps_{\mathbb{K}+2}$ to $\Lam=\veps_{\mathbb{K}+1}$.
Namely as in \cite{LMPS}, an exponential structure for indices emerges again.
Such a structure is needed to resolve the higher Mahlo operations $M_{k}$
\textit{inside} the lower classes $Mh_{i}$ for $i<k$.
A sequence $\la\alp_{k},\ldots,\alp_{N-1}\ra$ of ordinals $\alp_{i}<\Lam$
in this note
corresponds to a tower $t_{k}<\veps_{\mathbb{K}+2}$ of ordinals, where
$t_{i}=\Lam^{t_{i+1}}\alp_{i}$ and $t_{N-1}=\alp_{N-1}$.

\end{document}